\documentclass[12pt]{article}

\usepackage[margin=2cm]{geometry}
\usepackage{amsmath,amstext,amsfonts,amsthm,amssymb,bbm}
\usepackage{hyperref}

\def\mes{\operatorname{mes}}

\def\dist{\operatorname{dist}}

\def\ogamma{\overline{\gamma}}
\def\zen{\mathfrak z}

\theoremstyle{plain}
\newtheorem{thm}{Theorem}
\newtheorem{defin}{Definition}[section]

\newtheorem{prop}[defin]{Proposition}
\newtheorem{cl}[defin]{Claim}
\newtheorem{rmk}[defin]{Remark}

\newtheorem{lemma}[defin]{Lemma}
\newtheorem*{thm*}{Theorem}

\begin{document}

\author{Sasha Sodin\textsuperscript{1}}
\footnotetext[1]{Einstein Institute of Mathematics, Edmond J. Safra Campus, Givat Ram, The Hebrew University of Jerusalem,
Jerusalem, 9190401, Israel \& School of Mathematical Sciences,
Queen Mary University of London, 
Mile End Road, London E1 4NS, UK. Email: alexander.sodin@mail.huji.ac.il. Supported in part by a Royal Society Wolfson Research Merit Award (WM170012),  a Philip Leverhulme Prize of the Leverhulme Trust (PLP-2020-064), and the Or programme of the ISF under agreement number 4101/25.}
\title{On the Hausdorff measure of sets of non-Lyapunov behaviour, and a Jarn\'ik-type theorem for random Schr\"odinger operators.}
\maketitle 

\begin{abstract} We consider the growth of the norms of transfer matrices of  discrete Schr\"odinger operators with independent, identically distributed potential in one dimension. It is known that the set of energies at which the rate of exponential growth is slower than prescribed by the Lyapunov exponent is residual in the part of the spectrum at which the Lyapunov exponent is positive.  On the other hand, this exceptional set is of vanishing Hausdorff measure with respect to any gauge function $\rho(t)$ such that $\rho(t)/t$ is integrable at zero.

Here we show that this condition on $\rho(t)$ is optimal: for  potentials of sufficiently regular distribution, the set of energies at which the rate of exponential growth is arbitrarily slow has infinite Hausdorff measure with respect to any gauge function $\rho(t)$ such that $\rho(t)/t$ is non-increasing and not integrable at zero.

The main technical ingredient, possibly of independent interest, is a Jarn\'ik-type theorem describing the Hausdorff measure of the set of real numbers well-approximated by the eigenvalues of the Schr\"odinger operator. The proof of this result relies on the theory of Anderson localisation, and on the mass transference principle of Beresnevich--Velani.
\end{abstract}

\section{Introduction}

\paragraph{Sets of non--Lyapunov behaviour} Let $(v_k)_{k \geq 1}$ be an ergodic sequence, i.e.\ 
\[ v_k = v_k(\omega) = F(T^k \omega)~, \]
 where $(\Omega, \mathcal B, \mathbb P, T)$ is a probability space equipped with an ergodic transformation $T$, and $F: \Omega \to \mathbb R$ is a measurable function which we assume to satisfy $\mathbb E \log_+ |F(\omega)| < \infty$. We are interested in the Schr\"odinger operator $H$ defined by the random Jacobi matrix 
\begin{equation}\label{eq:defop}H = \left( \begin{array}{cccccc} 
v_1 & 1 & 0 &  & \\
1 & v_2 & 1 & 0 & \\
0 & 1 & v_3 & 1 & 0 & \\
 & 0 & 1 & v_4 & 1  & 0 \\
  &  &  & \ddots & \ddots & \ddots \\
\end{array} \right)~,\end{equation}
and in the associated transfer matrices 
\[ T_n(E) = \left( \begin{array}{cc} E - v_n & -1 \\ 1 & 0 \end{array}\right)~, \quad
\Phi_n(E) = T_n(E) \cdots T_2(E) T_1(E) \qquad (E \in \mathbb R)~. \]
Let 
\[ \gamma(E) = \lim_{n \to \infty} \frac1n \mathbb E \log \|\Phi_n(E) \|\]
be the Lyapunov exponent. A theorem of Furstenberg--Kesten \cite{FK} implies that for any given $E \in \mathbb R$ one has almost surely
\begin{equation}\label{eq:lyap}\lim_{n \to \infty} \frac1n   \log \|\Phi_n(E) \| = \gamma(E)~,\end{equation}
and consequently  the set $\Lambda$ of $E \in \mathbb R$ for which (\ref{eq:lyap}) fails is almost surely of zero measure. This leads to the question whether this exceptional set  $\Lambda$ is in fact empty. A positive result of Craig--Simon \cite{CS} asserts that almost surely 
\[ \forall E \in \mathbb R \,\,\, \limsup_{n \to \infty}  \frac1n   \log \|\Phi_n(E) \| \leq \gamma(E)~,\]
which rules out deviations in one direction (``up''). On the other hand,  deviations in the opposite direction (``down'') do occur when the Lyapunov exponent is positive. For $0 \leq \tau < 1$, let
\[ \Lambda^\tau = \left\{ E \in \mathbb R \, : \,   \, \liminf_{n \to \infty}  \frac1n   \log \|\Phi_n(E) \| \leq \tau \gamma(E)\right\}~, \]
and let $S$ be the essential spectrum of $H$.  In \cite{G75,G80a}, Goldsheid proved that for $\tau \geq 1/2$ the set $\Lambda^\tau$ is almost surely residual in $S \cap \{E\,:\, \gamma(E)>0\}$. Subsequently, he showed that this is in fact true for all $\tau \geq 0$ (see \cite[Theorem 2]{GS2}). Recently, a similar result was established by Gorodetski and Kleptsyn \cite{GorKl} for products of independent, identically distributed (i.i.d.) random matrices in SL$_2(\mathbb R)$, allowing for more complicated (non-linear) dependence of $T_n(E)$ on the parameter $E$. Thus the exceptional set $\Lambda$ is very thick in topological sense.

On the metric side, Gorodetski--Kleptsyn \cite{GorKl} proved that in the i.i.d.\  case the exceptional set is almost surely of zero Hausdorff  dimension. In the companion paper \cite{GS2} (devoted to the properties of exceptional sets for  operators  of the form (\ref{eq:defop}) as well as  their block matrix generalisations) we show, motivated by an idea of Simon \cite{Simon},  that the $\rho$-Hausdorff measure 
\[ \mes_\rho \Lambda = \lim_{\epsilon \to +0} \inf  \left\{ \sum_{j=1}^\infty \rho(\epsilon_j) \, : \, \Lambda \subset \bigcup_{j=1}^\infty (a_j - \epsilon_j, a_j + \epsilon_j)~, \, 0 < \epsilon_j \leq \epsilon \right\} \]
 vanishes for any gauge function $\rho(t)$ (i.e.\ a non-decreasing continuous function $\rho: [0, 1] \to \mathbb R_+$ with $\rho(0)=0$)  such that  
\begin{equation}\label{eq:intcond} \int_0^1 \frac{\rho(t)dt}{t} < \infty~.\end{equation}
The goal of the current paper is to demonstrate that the condition (\ref{eq:intcond}) is sharp. We assume:
\begin{equation}\label{eq:assum}
\begin{cases}
&\text{$v_k$ are independent, identically distributed random variables}\\
&\text{supported on a compact interval $J$, with  density bounded from above on $J$,}\\
& \text{and bounded from zero on any proper subinterval of $J$}~;
\end{cases}
\end{equation}
in this case, 
\[ S = [-2,2] + J = \bigcup_{E \in J} [E-2, E+2]~. \]
\begin{thm}\label{thm:1}
Assume (\ref{eq:assum}), and let $I \subset S$ be an interval. Then almost surely $\mes_\rho (\Lambda^\tau \cap I)  = \infty$ for any $\tau > 0$ and any gauge function $\rho(t)$ such that $\rho(t)/t$ is non-increasing and non-integrable on $(0, 1]$. 
\end{thm}

\begin{rmk} 
 It is likely that the statement is also true for $\tau = 0$. A proof along the lines of this paper would require a  stronger version of the results on Anderson localisation than what is currently available in the literature. 
\end{rmk}
\begin{rmk} While the assumptions (\ref{eq:assum}) on the distribution of $v_k$ are quite restrictive,  we do not believe that the result remains valid  for general ergodic sequences (say, with uniformly positive Lyapunov exponent). 
\end{rmk}
\begin{rmk}\label{rmk:kl-q-1} Kleptsyn and Quintino \cite{KlQ,Q} conjecture that the exceptional sets  are of positive logarithmic capacity.  Theorem~\ref{thm:1} provides plausible support for this conjecture, in view of the criterion of Frostman \cite{Fr}. 
\end{rmk}

\paragraph{Approximability by the eigenvalues of $H$}
To state the next results, we first summarise a few mostly known facts concerning the spectral properties of the random Schr\"odinger operator $H$. From now on, we assume that $v_k$ are independent, identically distributed and non-constant. In this case,  a theorem of Furstenberg \cite{F} ensures that $\gamma(E) > 0$ for any $E \in \mathbb R$; moreover, $\gamma(E)$ is continuous and tends to infinity at $\pm \infty$, whence $\inf_{E \in \mathbb R} \gamma(E) > 0$.
\begin{prop}\label{prop:loc} Assume that $v_k$ are independent, identically distributed, non-constant, and bounded. There exists $C > 0$ such that for any $\tau \in (0, 1)$ there exist almost surely
\[ (E_k, \psi_k) \in \mathbb R \times \mathbb \ell_2  \quad (k \in \mathbb N)  \]
and $K< \infty$ (all of them dependent on $\omega \in \Omega$ and, possibly, also on $\tau$) such that
\begin{enumerate}
\item[(L1)] $(\psi_k)$ form an orthonormal basis of $\ell_2$;
\item[(L2)] $H \psi_k = E_k \psi_k$;
\item[(L3)] $|\psi_k(x)| \leq \exp(- (1 - \tau) \gamma(E_k) | x - k|)$ for $x \in \mathbb N$ such that  $|x - k| \geq \max(\sqrt{k}, K)$;
\item [(L4)] for $\log k \leq k' < k$, $|E_k - E_{k'}| \geq \max(k, K)^{-C}$. 
\end{enumerate}
\end{prop}
That is, there is a basis of  eigenfunctions such that the $k$-th one is exponentially localised at the site $k$, and eigenvalues with close labels can not be too close to one another. For completeness, we sketch a proof in Section~\ref{s:props}.

\medskip
It turns out that the exceptional sets $\Lambda^\tau$ contain those $E$ which are well-approximable by the eigenvalues $E_k$. In the same  Section~\ref{s:props} we use Proposition~\ref{prop:loc} to prove:
\begin{prop}\label{prop:A}  Assume that $v_k$ are independent, identically distributed, non-constant, and bounded.  Let $I \subset S$ be a finite interval, and let $\ogamma_I = \max_{E \in I} \gamma(E)$. Then almost surely
\begin{equation}\label{eq:prop-A} \Lambda^{6 \tau} \cap I \supset \left\{ E \in  I \, : \, \# \left\{ k \geq 1 \, : \, |E - E_k| \leq e^{-2 \ogamma_I k} \right\} = \infty \right\}~.\end{equation}
\end{prop}
This leads us to the study of the metric properties of the set on the right-hand side of (\ref{eq:prop-A}). More generally, let $\alpha_k \searrow 0$, and let 
\[ \mathcal A_\alpha = \left\{ E \in \mathbb R \, : \, \# \left\{ k \geq 1 \, : \, |E - E_k| \leq  \alpha_k \right\} = \infty \right\}\]
be the set of $E$ well-approximable by the eigenvalues $E_k$.
We prove the following Jarn\'ik-type theorem which is perhaps of independent interest.
\begin{thm}\label{thm:J}
Assume (\ref{eq:assum}), and let $I \subset S$ be an interval, and $\alpha_k \searrow 0$.
\begin{enumerate}
\item If $\rho(t)$ is a gauge function such that $\sum_{k=1}^\infty \rho(\alpha_k) < \infty$, then almost surely  $\mes_\rho (\mathcal A_\alpha \cap I)= 0$;
\item If $\rho(t)$ is a gauge function such that $\rho(t)/t$ is non-increasing and $\sum_{k=1}^\infty \rho(\alpha_k) = \infty$, then almost surely $\mes_\rho (\mathcal A_\alpha \cap I)= \mes_\rho I$.
\end{enumerate}
\end{thm}
\begin{rmk} The classical Jarn\'ik theorem provides a similar-looking statement  for the case when $E_k$ are replaced with the rational numbers, numbered so that $p/q \prec p'/q'$ if $q < q'$.  See Beresnevich--Dickinson--Velani \cite{BDV} and references therein. 
\end{rmk} 
\begin{rmk}\label{rmk:kl-q-2} Quintino \cite{Q}  provides a sufficient condition for sets of the form appearing in the right-hand side of (\ref{eq:prop-A}) to be of positive logarithmic capacity. In view of Remark~\ref{rmk:kl-q-1} and Proposition~\ref{prop:A}, it would be interesting to check whether this condition is satisfied by the eigenvalues $E_k$.
\end{rmk}
\begin{proof}[Proof of Theorem~\ref{thm:1} (using Theorem~\ref{thm:J})] Let $\alpha_k = e^{-2 \ogamma k}$, so that
\[ \sum_k \rho(\alpha_k)\geq \int_1^\infty \rho(e^{-2 \ogamma s }) ds = \frac{1}{2 \ogamma} \int_0^{e^{-2 \ogamma}} \frac{\rho(t)dt}{t} = \infty~.  \]
By  the second item of Theorem~\ref{thm:J} almost surely $\mes_\rho  (\mathcal A_\alpha \cap I) = \mes_\rho I =  \infty$, whereas by Proposition~\ref{prop:A} almost surely $\Lambda^{6 \tau} \cap I \supset (\mathcal A_\alpha \cap I)$. Thus $\mes_\rho (\Lambda^{6 \tau} \cap I) = \mes_\rho I = \infty$ for any $\tau > 0$, as claimed.
\end{proof}

The mass transference principle of Beresnevich and Velani \cite{BV} allows to reduce the proof of Theorem~\ref{thm:J} to the special case of the Lebesgue measure, $\rho(t) = 2 t$ (a Khinchin-type theorem). We state this special case as a separate theorem.

\begin{thm}\label{thm:Kh} Assume (\ref{eq:assum}), and let $I \subset S$ be a finite interval, and $\alpha_k \searrow 0$. 
\begin{enumerate}
\item If $\sum \alpha_k < \infty$, then almost surely  $\mes  (\mathcal A_\alpha \cap I) = 0$;
\item If $\sum \alpha_k = \infty$, then almost surely $\mes  (\mathcal A_\alpha \cap I) = \mes I$.
\end{enumerate}
\end{thm}

\noindent Theorem~\ref{thm:Kh} is the main technical result of the current paper; its proof occupies Section~\ref{s:Kh} below.  

\begin{proof}[Proof of Theorem~\ref{thm:J} (using Theorem~\ref{thm:Kh})]
The first item is standard and does not rely on any properties of $E_k$. For any  $\epsilon, \delta > 0$, let $k(\epsilon)$ be such that $\alpha_{k(\epsilon)} \leq \epsilon$, and let $k'(\delta)$ be such that $\sum_{k \geq k'(\delta)} \rho(\alpha_k) \leq \delta$. Then 
\[ \mathcal A_\alpha \subset \bigcup_{k=\max(k(\epsilon), k'(\delta))}^\infty (E_k - \alpha_k, E_k + \alpha_k) \]
is a covering of $\mathcal A_\alpha$ by intervals of length $\leq \epsilon$ with 
\[ \sum_{k=\max(k(\epsilon), k'(\delta))}^\infty \rho(\alpha_k)< \delta~.\] 
Thus  $\mes_\rho  (\mathcal A_\alpha) \leq \delta$ for any $\delta > 0$, whence  $\mes_\rho  (\mathcal A_\alpha) =0$ and in particular $\mes_\rho  (\mathcal A_\alpha \cap I) = 0$. 
 
The second item follows from the second item of Theorem~\ref{thm:Kh} and the mass transference principle of Beresnevich and Velani  \cite[Theorem 2]{BV}.
\end{proof}

\section{Proof of Theorem~\ref{thm:Kh}}\label{s:Kh}

\paragraph{Preliminary reductions} The first item of  Theorem~\ref{thm:Kh} is a special case of the first item (proved above) of Theorem~\ref{thm:J}, therefore we turn to the second item of  Theorem~\ref{thm:Kh}.

Next, the following standard lemma allows us to assume that $\alpha_k \leq 1/k$, and we shall indeed make this assumption in the remainder of this section.
\begin{lemma} Let $\alpha_k \searrow 0$, $\sum_{k \geq 1} \alpha_k = \infty$. Then $\sum_{k \geq 1} \min(\alpha_k, 1/k) = \infty$.
\end{lemma}
\begin{proof}
We have:
\[\begin{split}
\sum_{2^m \leq k < 2^{m+1}} \min(\alpha_k, 1/k) 
&\geq \sum_{2^m \leq k < 2^{m+1}} \min(\alpha_k, 2^{-m-1}) 
\geq 2^m \min(\alpha_{2^{m+1}},  2^{-m-1}) \\
&= \frac12 \min(2^{m+1} \alpha_{2^{m+1}}, 1)~. \end{split}\]
If $2^m \alpha_{2^m} \geq 1$ for infinitely many values of $m$, then clearly 
\[ \sum_{k \geq 1} \min(\alpha_k, 1/k)  =  \sum_{m \geq 0} \sum_{2^m \leq k < 2^{m+1}}  \min(\alpha_k, 1/k)  \geq \frac12  \sum_{m \geq 1}\min(2^{m+1} \alpha_{2^{m+1}}, 1) = \infty~; \]
otherwise, 
\[ \sum_{m \geq 0} \min(2^{m+1} \alpha_{2^{m+1}}, 1) = \infty\]
by the Cauchy condensation test, whence again 
\[ \sum_{k \geq 1} \min(\alpha_k, 1/k) = \infty~. \qedhere \]
\end{proof}

\noindent Finally, recall that the integrated density of states of $H$ is defined as the  limit 
\[ \mathcal N(E) = \lim_{n \to \infty} \frac{1}{n} \, \mathbb E \# \left\{ \text{eigenvalues of $H_{[1,n]}$ in $(-\infty, E]$} \right\}~,\]
where $H_{[1,n]}$ is the restriction of $H$ to $[1, n]$ (the submatrix  formed by the first $n$ rows and columns of (\ref{eq:defop})). The Wegner estimate \cite{Weg} (presented in textbook form in \cite{AW}) asserts that if the distribution of $v_1$ is absolutely continuous with density bounded by $A$, then so is $\mathcal N$:  
\[ \mathcal N(E) = \int_{-\infty}^E \mathcal N'(E') dE'~, \quad \mathcal N' \leq A~. \]
Without loss of generality, we may assume that $A \geq 1$.

Further, a lower Wegner estimate \cite{Weg,Jeske,HM} asserts that under the assumptions (\ref{eq:assum})  $\mathcal N'$ is bounded from $0$ on any interval  $I$ contained in the interior of $S$. It is sufficient to prove the theorem for such intervals $I$; thus in the sequel we assume that $I$ is contained in the interior of $S$ and we denote $a_I = \inf_{E \in I} \mathcal N'(E)$.

\paragraph{Reduction to two main claims}
\leavevmode 
% [removed, formerly red:] Let $H_m$ be the restriction of $H$ to the interval $[4^m, 2 \times 4^m)$, i.e.\ the submatrix of (\ref{eq:defop}) formed by the rows and columns with indices $4^m \leq k < 2\times 4^m$, and let $\sigma(H_m)$ be the spectrum of $H_m$.
Let 
\[ L_m = 2^{m + \lceil m/2 \rceil}~, \quad n_m = 4^m / L_m = 2^{\lfloor m/2 \rfloor}~.\]
Divide the interval $[4^m, 2 \times 4^m)$ into $n_m$ consecutive blocks $\Lambda_{m,1}, \cdots, \Lambda_{m, n_m}$ of length $L_m$, let $H_m^{(j)}$ be the restriction of $H$ to $\Lambda_{m,j}$, i.e.\ the submatrix of (\ref{eq:defop}) formed by the rows and columns with indices in $\Lambda_{m,j}$, and let 
\[ H_m = \bigoplus_{j=1}^{n_m} H_m^{(j)}~, \quad \sigma(H_m) = \bigcup_{j=1}^{n_m} \sigma(H_m^{(j)})~.\]
That is, $H_m$ is obtained from the restriction of $H$ to $[4^m, 2 \times 4^m)$ by removing the $n_m - 1$ bonds connecting adjacent blocks. The operators $H_m^{(j)}$ are jointly independent (for all $m$ and $j$); the decomposition into blocks is used in the proof of Claim~\ref{cl:induc} below. Denote 
\[ B_I = \left\{ E \in I \, : \, \# \left\{ k \, : \, |E - E_k| < \alpha_k \right\} < \infty \right\} \]
and 
\[B_I' = \bigcup_{m \geq 1} B_{m,I}'~, \quad \text{where} \quad
B_{m,I}' = \left\{ E \in I \, : \, \forall m' \geq m \, : \, \dist(E, \sigma(H_{m'})) > \frac12  \alpha_{2 \times 4^{m'}}  \right\}~.\]
We shall prove:
\begin{cl}\label{cl:decoup} Almost surely $\mes (B_I \setminus B_I') =0$.
\end{cl}
That is, the set $B_I$ of numbers poorly approximable by the eigenvalues of $H$ is almost contained in the set $B_I'$ of numbers poorly approximable by the eigenvalues of the chunks $H_m$.
The advantage of $B_I'$ over $B_I$ is that the spectra $\sigma(H_m)$ appearing in the definition of the former one are jointly independent, whence for any given $E \in I$ the event $\{ E \in B_I' \}$ lies in the tail sigma-algebra and thus satisfies the zero-one law. Consider the deterministic set 
\[ B_I '' = \left\{ E \in I \, : \, \mathbb P\{ E \in B_I'\} = 1\right\}~; \]
then almost surely $\mes (B_I' \triangle B_I'') =0$. Thus we need to show that $\mes B_I '' = 0$.
We shall prove that typically either the measure of $B_{m+1}'$ is already small, or  the measure of $B_{m}' \subset B_{m+1}'$ is significantly smaller than that of $B_{m+1}'$, i.e.\ a sizeable fraction of the numbers poorly approximable by $\sigma(H_{m'})$ with indices $m' > m$ is well approximable by $\sigma(H_m)$. Formally,
\begin{cl}\label{cl:induc} There exists $c_0(a_I, A)> 0$ (not dependent on $\mes I$!) such that 
for any $0 < \zeta \leq c_0(a_I, A)$, and any $m \geq m_0$ (a deterministic number possibly depending on all the parameters, including $\zeta$), 
\[ \begin{split} &\mathbb P \left\{ \mes B'_{m+1,I} \geq (1-\zeta) \mes I~, \quad
\mes (B'_{m+1,I} \setminus B'_{m,I}) \leq  \zeta \mes I \, 4^m \alpha_{2 \times4^m}  \right\} \\
&\qquad\qquad\leq  \exp\left(- c_0(a_I, A) (\mes I)^2 \, 2^{m/2}\right)~. \end{split}\]
% [removed, formerly red:] e^{- c_0(a_I, A) \mes I \, 4^m}
\end{cl}
The claim implies that almost surely we have $\mes B_{m, I}' < (1-\zeta) \mes I$ for all $m$. Indeed, if this inequality would fail for a certain $m = m_1$, it would definitely fail for all $m \geq m_1$. By the first Borel-Cantelli lemma (the probabilities in Claim~\ref{cl:induc} are summable in $m$), this would imply that 
\[ \mes (B'_{m+1,I} \setminus B'_{m,I}) >  \zeta \mes I \, 4^m \alpha_{2 \times 4^m}\]
for all sufficiently large $m$. But then an application of the Cauchy condensation test would yield
\[ \mes B'_{m, I} \leq \mes I - \sum_{m=m_1}^\infty  \zeta \mes I \,  4^m \alpha_{2 \times 4^m} = -\infty~,\]
which can not be true. Thus indeed $\mes B_{m, I}' < (1-\zeta) \mes I$ for all $m$, and by similar reasoning also  $\mes B_{m, \tilde I}' < (1-\zeta) \mes \tilde I$ for any (deterministic)  subinterval $\tilde I \subset I$.

Now we can prove that $\mes B_I'' = 0$.  If this were not the case, $B_I''$ would have a Lebesgue point $E^*$, and then we would be able to find a (deterministic) neighbourhood  $\tilde I \subset I$ of $E^*$ such that $\mes B_{\tilde I}''  > (1-\zeta) \mes \tilde I$. But then almost surely we would have $\mes B'_{\tilde I} > (1-\zeta) \mes \tilde I$, whence $\mes B'_{m, \tilde I} > (1-\zeta) \mes \tilde I$ for sufficiently large $m$, which contradicts what we have proved. Thus the theorem is reduced to Claims~\ref{cl:decoup} and \ref{cl:induc}, stated above and proved below.

\paragraph{Proof of Claim~\ref{cl:decoup}}

\leavevmode Let $G_m$ be the set of those $k \in [4^m, 2 \times 4^m)$ which lie at distance at least $2^{m+1}$ from both endpoints of the block $\Lambda_{m,j}$ containing $k$; thus
\begin{equation}\label{eq:countgood} \# \left( [4^m, 2 \times 4^m) \setminus G_m \right) \leq 2^{m+2} n_m \leq 4 \times 2^{3m/2}~.\end{equation}
First, observe that almost surely for $m$ large enough (depending on the realisation)
\begin{equation}\label{eq:comp-sp} \max_{ k \in G_m} \dist(E_k, \sigma(H_m)) \leq e^{-c 2^m} < \frac12 \times 8^{-m}~. \end{equation}
% [removed, formerly red:] 4^m + 2^m \leq k < 2\times4^m - 2^m
Indeed,  for  $k \in G_m$ let $\psi_k^m$ be the restriction of $\psi_k$ (the eigenfunction of $H$ corresponding to the eigenvalue $E_k$)  to  the block $\Lambda_{m,j}$ containing $k$. Every site $x \notin \Lambda_{m,j}$ satisfies $|x - k| \geq 2^{m+1} > \sqrt{k}$. Then by item (L3) of Proposition~\ref{prop:loc} 
% [removed, formerly red:] $4^m + 2^m \leq k < 2\times4^m - 2^m$
% [removed, formerly red:] $[4^m, 2 \times 4^m)$
\[ \| \psi_{k}^m\| \geq 1 - \exp(-c 2^m)~, \quad \| ( H_m^{(j)} - E_k) \psi_k^m \| \leq  \exp(-c 2^m)~,\]
% [removed, formerly red:] H_m
which implies (\ref{eq:comp-sp}) since $\sigma(H_m^{(j)}) \subset \sigma(H_m)$. 

Let 
\[ \sigma^g(H_m) = \left\{ E \in \sigma(H_m) \, : \, \exists \,  k \in G_m : |E - E_k| \leq  e^{-c 2^m} \right\}~, \]
% [removed, formerly red:] 4^m + 2^m \leq k < 2\times4^m - 2^m
and let $\sigma^b(H_m) = \sigma(H_m) \setminus \sigma^g(H_m)$.
By item (L4) of Proposition 1.4, for any $E \in \sigma^g(H_m)$ the index $k$ appearing in the definition of $\sigma^g(H_m)$ is unique; hence  $\# \sigma^g(H_m) \geq \# G_m$, and, by (\ref{eq:countgood}), $\# \sigma^b (H_m) \leq 4 \times 2^{3m/2}$.  
% [removed, formerly red:] $\# \sigma^g(H_m) \geq 4^m - 2^{m+1}$ and $\# \sigma^b (H_m) \leq 2^{m+1}$

Denote:
\[ \begin{split}
A_I^1 &= \left\{ E \in I \, : \, \# \left\{ m \, : \, \alpha_{2 \times 4^m} \leq   8^{-m}~, \,\,
\dist(E, \sigma(H_m)) \leq \frac12 \alpha_{2 \times 4^m} \right\} = \infty \right\}~,\\
A_I^2 &= \left\{ E \in I \, : \, \# \left\{ m \, : \, \alpha_{2 \times 4^m} >   8^{-m}~, \,\,
 \min_{k \in G_m}  |E-E_k| <  \alpha_{2 \times 4^m} \right\} = \infty \right\}~,\\
A_I^3 &= \left\{ E \in I \, : \, \# \left\{ m \, : \,  \dist(E, \sigma^b(H_m)) \leq 
\frac12 \alpha_{2 \times 4^m} \right\} = \infty \right\}~.
\end{split}\]
% [removed, formerly red:] (In the original text, the minimum in the definition of $A_I^2$ was taken over $4^m + 2^m \leq k < 2 \times 4^m - 2^m$.)
Then almost surely $I \setminus B_I' \subset A_I^1 \cup A_I^2 \cup A_I^3$. 
We have: $A_I^2 \cap B_I = \varnothing$, while the convergence of the series
\[\begin{aligned}
\sum_{m \, : \, \alpha_{2 \times 4^m} \leq 8^{-m}} \# \sigma(H_m) \times \alpha_{2 \times 4^m} 
&\leq \sum_{m} 4^m   \times 8 ^{-m}  &< \infty \\
\sum_m \# \sigma^b(H_m) \times \alpha_{2\times 4^m} &\leq \sum_m  4 \times 2^{3m/2} \times \frac{1}{2 \times 4^m} &< \infty 
% [removed, formerly red:] 2^{m+1}
\end{aligned}\]
(where in the second display we have used that $\alpha_k \leq 1/k$) implies, with the help of  the first Borel--Cantelli lemma, that $\mes A_I^1 = \mes A_I^3 = 0$. \qed

\paragraph{Proof of Claim~\ref{cl:induc}}
The random sets $\sigma(H_m)$ and $B'_{m+1,I}$ are independent, hence it suffices to show that for any (given) $B \subset I$ of measure
$\mes B \geq (1- \zeta) \mes I$
\begin{equation}\label{eq:need-ind} \mathbb P \left\{ \mes( B  \, \cap \, \Delta_{m, I} ) \leq  \zeta \mes I \,  4^m \alpha_{2\times4^m} \right\}  \leq  \exp\left(- c_0(a_I, A) (\mes I)^2 \, 2^{m/2}\right)~,\end{equation}
% [removed, formerly red:] \, \setminus \,
% [removed, formerly red:] e^{-c_0(a_I,A) \mes I \, 4^m}
where
\[ \Delta_{m,I} = \left\{ E \in I \, : \, \dist(E, \sigma(H_m)) <\frac12\alpha_{2 \times 4^m} \right\}~.\]
(Once this is proved, we can take $B = B'_{m+1, I}$ while conditioning on this set; note that $B'_{m+1,I} \setminus B'_{m,I} \supset B'_{m+1,I} \cap \Delta_{m,I}$.)
We need  the elementary
\begin{lemma}\label{l:cov} Let $I$ be an interval, and let  $B \subset I$ be a measurable set. Then for any $\theta > 0$
\[ \mes A^\theta \leq 4 (\mes(I \setminus B) + \theta)~, \quad \text{where} \quad A^\theta = \left\{ E \in I \, : \, \mes ((E - \theta, E+\theta) \cap B) \leq \theta \right\}~.\]
\end{lemma}
\begin{proof}
Let
\[ a_1 = \inf \{ a \in A^\theta\}~, \quad a_2 = \inf \{ a \in A^\theta~, \, a \geq a_1 + \theta\}~, \quad
a_3 = \inf \{a \in A^\theta~, \, a \geq a_2 + \theta\}~,\cdots\]
Then 
\[\begin{split}
\mes A^\theta &\leq \mes \bigcup_j (a_j - \theta, a_j + \theta)
\leq \sum_j \mes (a_j - \theta, a_j + \theta) \\
&\leq 2 \sum_j \mes ((a_j - \theta, a_j + \theta) \setminus B) \leq 4 \mes(I \setminus B) + 4 \theta ~.\qedhere
\end{split}\]
\end{proof}
Now we proceed with the proof of the claim. Applying Lemma~\ref{l:cov} to our set $B$ with $\theta = \frac12 \alpha_{2 \times 4^m}$, we have (for sufficiently large $m$) $\mes A^\theta \leq 5 \zeta \mes I$. 

Let 
\begin{equation}\label{eq:defK} K = \left\lceil \frac{32 A^2}{a_I} \right\rceil~, \quad \text{and assume that} \quad 0 < \zeta \leq \min \left( \frac{a_I}{40 A}, \, \frac{a_I}{8K} \right)~. \end{equation}
Note that $K \geq 32 A \geq 32$, since $a_I \leq A$ and $A \geq 1$. Divide $I$ into $\lceil K \, 4^{m}\mes I \rceil$ intervals $J$ of length $\leq 4^{-m}/K$ so that all the intervals except for the last one are exactly of length $4^{-m}/K$, and let 
\[ Z = \# \left\{ J \, : \, J \cap (\sigma(H_m) \setminus A^\theta) \neq \varnothing \right\}~.\]
In each of the $Z$ intervals $J$ counted here, choose a point $E_J \in J \cap (\sigma(H_m) \setminus A^\theta)$. By the definitions of $A^\theta$ and of $\Delta_{m,I}$, 
\[ \mes \left((E_J - \theta, E_J + \theta) \cap B \right) > \theta~, \quad (E_J - \theta, E_J + \theta) \cap I \subset \Delta_{m,I}~.\]
As $2 \theta = \alpha_{2 \times 4^m} \leq \frac12 4^{-m}$, an interval of length $2\theta$ intersects at most $\lceil K/2 \rceil + 1 \leq K$ of the intervals $J$, hence each point is covered by at most $K$ of the neighbourhoods $(E_J - \theta, E_J + \theta)$. Therefore
\begin{equation}\label{eq:Zmes}
\mes (B \cap \Delta_{m,I}) \geq \frac{1}{K} \sum_J \mes \left((E_J - \theta, E_J + \theta) \cap B \right) \geq \frac{\theta Z}{K} = \frac{\alpha_{2 \times 4^m}}{2K} Z~.
\end{equation}
In view of (\ref{eq:defK}) and (\ref{eq:Zmes}), the estimate (\ref{eq:need-ind}) will follow once we show that for sufficiently large $m$
\begin{equation}\label{eq:Zprob}
\mathbb P \left\{ Z \leq \frac{a_I}{4} \mes I \, 4^m \right\} \leq \exp \left( - \frac{a_I^2}{8} (\mes I)^2 \, 2^{\lfloor m/2 \rfloor} \right)~;
\end{equation}
indeed, on the complement of this event $\mes (B \cap \Delta_{m,I}) > \frac{a_I}{8K} \mes I \, 4^m \alpha_{2 \times 4^m} \geq \zeta \mes I \, 4^m \alpha_{2 \times 4^m}$. Thus Claim~\ref{cl:induc} holds with $c_0(a_I, A) = \min(\frac{a_I}{40A}, \frac{a_I}{8K}, \frac{a_I^2}{16})$. We prove (\ref{eq:Zprob}) in two steps: first we bound the expectation of $Z$ from below, and then we show that $Z$ is concentrated about its expectation.

\smallskip\noindent
\textit{The expectation of $Z$.} For a Borel set $S \subset \mathbb R$, let 
\[ N_S^{(j)} = \# ( \sigma(H_m^{(j)}) \cap S)~, \quad N_S = \sum_{j=1}^{n_m} N_S^{(j)} = \# ( \sigma(H_m) \cap S) \]
(the eigenvalues are counted with multiplicity). The number of intervals $J$ containing a point of $\sigma(H_m)$ is equal to $\sum_J (N_J - (N_J - 1)_+) = N_I - \sum_J (N_J - 1)_+$, and each interval $J$ which contains a point of $\sigma(H_m)$ but no point of $\sigma(H_m) \setminus A^\theta$ contains at least one point of $\sigma(H_m) \cap A^\theta$. Therefore
\begin{equation}\label{eq:Zlower}
Z \geq N_I - \sum_J (N_J - 1)_+ - N_{A^\theta}~.
\end{equation}
We estimate the expectations of the three terms.

First, choose a large natural $M$, and compare the restriction $H_{[1, M L_m]}$ of $H$ to $[1, M L_m]$ with the direct sum of its restrictions to $M$ consecutive intervals of length $L_m$, which are independent copies of $H_m^{(1)}$. By interlacing,
\[ \mathbb E \, \# ( \sigma(H_{[1, M L_m]}) \cap I) \leq M \left(\mathbb E N_I^{(1)} + 2 \right)~,\]
whence, dividing by $M L_m$, letting $M \to \infty$, and using that $\mathcal N'  \geq a_I$ on $I$, we obtain $\mathbb E N_I^{(1)} \geq a_I \mes I \, L_m - 2$. Thus
\begin{equation}\label{eq:EZ1}
\mathbb E N_I \geq a_I \mes I \, 4^m - 2 n_m~.
\end{equation}

Second, let $J$ be one of the intervals of the partition. By the Minami bound (\cite{Min}; see Combes--Germinet--Klein \cite{CGK} for the form we use here), applied to the {interval} $J$ and to the operator $H_m^{(j)}$,
\[ \mathbb P \left\{ N_J^{(j)} \geq r \right\} \leq \frac{Q^r}{r!}~, \quad \text{where} \quad Q = A \mes J \, L_m \leq \frac{A L_m}{K 4^m} \leq 1~,\]
whence 
\[ \mathbb E N_J^{(j)} (N_J^{(j)} - 1) = \sum_{r \geq 2} 2 (r-1)  \mathbb P \left\{ N_J^{(j)} \geq r \right\} \leq Q^2 \sum_{r \geq 2} \frac{Q^{r-2}}{(r-2)!} \leq e \, Q^2~.\]
The Wegner estimate yields $\mathbb E N_J^{(j)} \leq Q$. As the blocks are independent, 
\[ \mathbb E N_J (N_J - 1) = \sum_j \mathbb E N_J^{(j)} (N_J^{(j)} - 1) + \sum_{j \neq j'} \mathbb E N_J^{(j)} \, \mathbb E N_J^{(j')} 
\leq e \, n_m Q^2 + n_m^2 Q^2 \leq \frac{4 A^2}{K^2}~.\]
Using the inequality $(N-1)_+ \leq \frac12 N(N-1)$, summing over the $\lceil K \, 4^m \mes I \rceil \leq 2 K \, 4^m \mes I$ intervals $J$, and recalling (\ref{eq:defK}), we obtain
\begin{equation}\label{eq:EZ2}
\mathbb E \sum_J (N_J - 1)_+ \leq \frac{4 A^2}{K} \mes I \, 4^m \leq \frac{a_I}{8} \mes I \, 4^m~.
\end{equation}

Third, the Wegner estimate \cite{Weg} is valid for an arbitrary Borel set (see \cite{AW}): $\mathbb E N_S^{(j)} \leq A \mes S \, L_m$. Applying it to $S = A^\theta$ and recalling that $\mes A^\theta \leq 5 \zeta \mes I$, we obtain using (\ref{eq:defK})
\begin{equation}\label{eq:EZ3}
\mathbb E N_{A^\theta} \leq A \mes A^\theta \, 4^m \leq 5 A \zeta \mes I \, 4^m \leq \frac{a_I}{8} \mes I \, 4^m~.
\end{equation}
Note that, unlike the Minami bound, the Wegner estimate does not require the set to be an interval. 

Combining (\ref{eq:Zlower}) with (\ref{eq:EZ1}), (\ref{eq:EZ2}) and (\ref{eq:EZ3}), we obtain for sufficiently large $m$
\begin{equation}\label{eq:EZ}
\mathbb E Z \geq \frac{3 a_I}{4} \mes I \, 4^m - 2 n_m \geq \frac{a_I}{2} \mes I \, 4^m~.
\end{equation}

\smallskip\noindent
\textit{Concentration.} The random variable $Z$ is a function of the $n_m$ independent blocks of potentials $(v_k)_{k \in \Lambda_{m,j}}$, $1 \leq j \leq n_m$. If the potentials in one block $\Lambda_{m,j}$ are changed, then $\sigma(H_m)$ changes only by replacing the $L_m$ eigenvalues of $H_m^{(j)}$ with $L_m$ other numbers, whence $Z$ changes by at most $L_m$. By the bounded differences inequality (McDiarmid \cite{McD}), for any $t > 0$
\[ \mathbb P \left\{ Z \leq \mathbb E Z - t \right\} \leq \exp \left( - \frac{2 t^2}{n_m L_m^2} \right) = \exp \left( - \frac{2 t^2}{4^m L_m} \right)~.\]
Taking $t = \frac{a_I}{4} \mes I \, 4^m$ and using (\ref{eq:EZ}) and $4^m / L_m = 2^{\lfloor m/2 \rfloor}$, we obtain (\ref{eq:Zprob}). This concludes the proof of Claim~\ref{cl:induc} and of Theorem~\ref{thm:Kh}. \qed

\section{Localisation and the exceptional sets}\label{s:props}

\begin{proof}[Sketch of proof of Proposition~\ref{prop:loc}]
Let $\tau > 0$. It is well-known that there exist  (almost surely)
\[ (E_k, \psi_k, \zen_k) \in \mathbb R \times \mathbb \ell_2 \times \mathbb N  \quad (k \in \mathbb N) 
 \]
and $K < \infty$ such that  
\begin{enumerate}
\item[(L1)] $(\psi_k)$ form an orthonormal basis of $\ell_2$;
\item[(L2)] $H \psi_k = E_k \psi_k$;
\item[(L$\tilde3$)] $|\psi_k(x)| \leq \exp(- (1 - \tau/2) \gamma(E_k) | x - \zen_k|)$ for $|x - \zen_k| \geq \max(\sqrt{\zen_k} / 10 , K)$.
\end{enumerate}
The first two items are explicitly stated, for example, in \cite{GorKl}; the third one follows from the first two using the condition on resonant boxes as established by multi-scale analysis in \cite{CKM} or by a single-scale argument in \cite{JZ}. 

To upgrade (L$\tilde 3$) to (L3), it remains to show that one can reorder the eigenfunctions so that $\zen_k = k$. To this end, it is sufficient to verify that for large $L$
\begin{equation}\label{eq:drjls+}  L - \sqrt{L} / 5 \leq \# \left\{ k \in \mathbb N \, : \, \zen_k \leq L \right\} \leq L + \sqrt{L} / 5~:\end{equation}
having (\ref{eq:drjls+}), we reorder $\psi_k$ so that $\zen_1 \leq \zen_2 \leq \zen_3 \leq \cdots$. The proof of (\ref{eq:drjls+}) follows the argument of del Rio, Jitomirskaya, Last, and Simon \cite[Theorem 7.1]{dRJLS}, replacing the input (SULE) assumed in \cite{dRJLS} with the stronger one given by (L$\tilde3$) above. We omit the details.

The proof of the lower bound (L4) on the spacings between the eigenvalues requires  a lemma that is also used in the proof of Proposition~\ref{prop:A} below.
\begin{lemma}\label{l:cot} In the setting of Proposition~\ref{prop:loc}, let $H_{[a, 2k + b-1]}$ be the restrictions of $H$ to the intervals $[a, 2k + b-1]$, where $a, b = 1,2$. Almost surely there exists $k_0 < \infty$ such that for $k \geq k_0$ each of these four  operators has an  eigenvalue $E(a,b)$ with eigenvector $\psi_k(a, b)$ satisfying
\begin{equation}\label{eq:closeev1}
|E_k(a,b) - E_k| \leq \exp(-2(1 - 2\tau)\gamma(E_k)k)~, \, \| \psi_k(a, b) - \psi_k|_{[a, 2k + b-1]}\| \leq \exp(- (1 - 2\tau) \gamma(E_k)k)~.
\end{equation} 
\end{lemma}
The precise estimates (\ref{eq:closeev1})  will be important in the proof of Proposition~\ref{prop:A}; here (in the proof of Proposition~\ref{prop:loc}) it is sufficient for us that the distances between $E_k(1,1), E_{k'}(1,1)$ and $E_k, E_{k'}$, respectively, decay faster than any power of $k$, and that $\psi_k(1, 1)$ and $\psi_{k'}(1,1)$ can not coincide since they are close to two orthogonal functions. According to a result of Bourgain \cite{Bourg}, almost surely the spacings between the eigenvalues of  $H_{[1, 2k]}$ are $\geq k^{-C}$, when $k$ is large enough (the result of Bourgain is stated for Bernoulli distribution, but the argument applies in the general case; see \cite{Cot}). Thus also 
\[ |E_k - E_{k'}| \geq k^{-C} - 2\exp(- c k) \geq k^{-C-1}~. \]
This concludes the proof of (L4).
\end{proof}

\noindent To complete the proof of Proposition~\ref{prop:loc}, it remains to prove Lemma~\ref{l:cot}.

\begin{proof}[Proof of Lemma~\ref{l:cot}]
This proof relies significantly on the arguments from the work of Haeming \cite{Cot}, to which we also refer for a detailed exposition. First observe (using the already proved item (L3) of Proposition~\ref{prop:loc})  that the trimmed eigenfunctions $\tilde \psi_k = \psi_k|_{[2, 2k-2]}$ satisfy 
\[ \| (H_{[a, 2k + b-1]} - E_k) \tilde \psi_k \| \leq C_1 \exp(- (1- \tau) \gamma(E_k)  k)~. \] 
This estimate ensures that the distance from $E_k$ to the spectrum of $H_{[a, 2k + b-1]}$ is at most 
 \[ \dist(E_k, \sigma(H_{[a, 2k + b-1]})) \leq C_2 \exp(- (1 - \tau) \gamma(E_k) k)~; \]
note that the exponent is almost twice smaller than what we claimed in (\ref{eq:closeev1}).  

To upgrade the proved bound to (\ref{eq:closeev1}), we first appeal  (as in the proof of Proposition~\ref{prop:loc} above) to a result of Bourgain \cite{Bourg} asserting  that the spacings between the eigenvalues of $H_{[a, 2k + b-1]}$ are lower-bounded by $k^{-C}$ (for $k$ large enough, depending on the realisation). This implies that the eigenvectors $\psi(a, b)$ that we have constructed also satisfy
\begin{equation}\label{eq:pert-psi} \| \psi(a, b) - \tilde \psi_k \| \leq C_2 k^C \exp(- (1 - \tau) \gamma(E_k) k) \leq  \frac12 \exp(- (1 - 2\tau) \gamma(E_k) k)~,\end{equation}
which proves the second part of (\ref{eq:closeev1}).

It remains to establish the first part of (\ref{eq:closeev1}), and to this end  we use first-order perturbation theory. Choose $N = 100 k$ to ensure that $H_{[1, N]}$ has an eigenvalue at distance $e^{-30 \gamma(E_k) k}$ from $E_k$. Construct the operators $H(t; a, b)$, $0 \leq t \leq 1$, interpolating linarly between  $H(0; a, b) = H_{[1,N]}$ and 
\[ H(1; a, b) = H_{[1, a-1]} \oplus H_{[a, 2k+b-1]} \oplus H_{[2k+b,N]}~. \]
Repeating the arguments presented above for $H_{[a,2k+b-1]}$ (which corresponds to $t=1$), we show that  for each $0 \leq t \leq 1$  the operator $H(t; a, b)$ has an eigenfunction $\psi_k(t;a,b)$ satisfying the counterpart of the bound (\ref{eq:pert-psi}) and an eigenvalue $E_k(t, a,b)$ satisfying 
\[ |E_k(t; a,b) - E_k| \leq C_2 \exp(-(1-\tau)\gamma(E_k)k)~.\]
Denoting by $\psi_k(t; a, b; \cdot)$ the components of $\psi_k(t; a, b)$, and applying the Feynman--Hellmann formula, we have:
\[\begin{split} 
&|E_k(1; a,b) - E_k(0; a,b)| \\
&\quad \leq \max_{0 \leq t \leq 1} |\langle \left( \frac{d}{dt} H(t; a,b) \right) \psi_k(t; a,b), \psi_k(t, a,b) \rangle| \\
&\quad \leq  \max_{0 \leq t \leq 1}  \big[ |\psi_k(t; a, b; a-1)| \times |\psi_k(t; a, b; a)| + |\psi_k(t; a, b; 2k+b-1)| \times |\psi_k(t; a, b; 2k+b)| \big] \\
&\quad \leq \exp(- 2 (1 - 2\tau) \gamma(E_k) k)~,
\end{split}\]
as claimed.
\end{proof}

\noindent Having established Lemma~\ref{l:cot}, we can turn to 
\begin{proof}[Proof of Proposition~\ref{prop:A}]
It suffices to prove that for $k$ large enough (depending on the realisation)
\begin{equation}\label{eq:need-A}
|E - E_k| \leq \exp(-2 \ogamma_I k) \Longrightarrow \| \Phi_{2k}(E)\| \leq e^{12\tau \gamma(E_k) k}~.
\end{equation}
We can assume that the interval $I$ containing the numbers $E, E_k$ is sufficiently short to ensure that 
\begin{equation}\label{eq:vargamma}\max_{\tilde E \in I} \gamma(\tilde E)  \leq (1 +\tau)  \min_{\tilde E \in I} \gamma(\tilde E)~. \end{equation}
Set $n = 2k$. The matrix entries of $\Phi_n(E)$ are exactly the determinants of the  four matrices $H_{[a, n + b-1]}$, $a,b=1,2$; Lemma~\ref{l:cot} shows that each of these has a zero within distance $e^{-(1 - 2\tau)\gamma(E_k) n}$ from $E_k$.
Now, a locally uniform version of the result of Craig--Simon \cite{CS} discussed in the introduction (see e.g.\ \cite{GS2}) implies that for $n$ large enough 
\[ \forall E \in  I \, : \,  \| \Phi_n(E) \| \leq \exp((1+\tau)\gamma(E) n)  \leq \exp((1 +2\tau)\gamma(E_k) n) ~,\]
whence by Markov's polynomial inequality
\[ \| \Phi_n'(E) \| \leq C_I n^2 \exp((1 +2\tau)\gamma(E_k) n) \leq \exp((1 +3\tau)\gamma(E_k) n) ~.\]
Thus for  any $E$ such that 
\[ |E - E_k| \leq \exp(-2 \ogamma_I k) = \exp(- \ogamma_I n)\]
we have:
\[ \| \Phi_n(E)\| \leq (e^{- \ogamma_I n}+ e^{-(1 - 2\tau)\gamma(E_k)n}) \exp((1 +3\tau)\gamma(E_k)n)  \leq \exp(6 \tau n)~,\]
as claimed in (\ref{eq:need-A}). 
\end{proof}

\paragraph{Acknowledgement} I am grateful to Anton Gorodetski for helpful comments on a preliminary version of this note, and in particular for bringing the work \cite{Q} to my attention. A mistake in the application of the Minami estimate has been corrected with the help of the Claude Team plan for Scientists.


\begin{thebibliography}{cc}
\bibitem{AW} Aizenman, Michael; Warzel, Simone.
 Random operators. Disorder effects on quantum spectra and dynamics. Graduate Studies in Mathematics, 168. American Mathematical Society, Providence, RI, 2015. xiv+326 pp. 
\bibitem{BV}  Beresnevich, Victor; Velani, Sanju.
A mass transference principle and the Duffin-Schaeffer conjecture for Hausdorff measures. 
Ann.\ of Math.\ (2) 164 (2006), no.~3, 971--992. 
\bibitem{BDV} Beresnevich, Victor; Dickinson, Detta; Velani, Sanju.
Measure theoretic laws for lim sup sets.  
Mem.\ Amer.\ Math.\ Soc.~179 (2006), no.~846, x+91 pp. 
\bibitem{Bourg} Bourgain, Jean.
 On eigenvalue spacings for the 1-D Anderson model with singular site distribution. 
Geometric aspects of functional analysis, 71--83, Lecture Notes in Math., 2116, Springer, Cham, 2014. 
\bibitem{CKM} Carmona, Ren\'e; Klein, Abel; Martinelli, Fabio.
 Anderson localization for Bernoulli and other singular potentials. Comm. Math. Phys. 108 (1987), no. 1, 41-66.
\bibitem{CGK} Combes, Jean-Michel; Germinet, Fran{\c c}ois; Klein, Abel.
Generalized eigenvalue-counting estimates for the Anderson model.  
J.\ Stat.\ Phys.\ 135 (2009), no.~2, 201--216.

\bibitem{CS}  
Craig, W.; Simon, B. 
Subharmonicity of the Lyaponov index. 
Duke Math.\ J.~50 (1983), no.~2, 551--560.
\bibitem{dRJLS} del Rio, R.; Jitomirskaya, S.; Last, Y.; Simon, B.
Operators with singular continuous spectrum. IV. Hausdorff dimensions, rank one perturbations, and localization. 
J.\ Anal.\ Math.~69 (1996), 153--200.
\bibitem{Fr} Frostman, O.
Potentiel d'\'equilibre et capacit\'e des ensembles. Avec quelques applications
a la th\'eorie des fonctions. PhD Thesis, Lund (1935). 
\bibitem{F} Furstenberg, H.
 Noncommuting random products. 
Trans.\ Amer.\ Math.\ Soc.\ 108 (1963), 377--428.
\bibitem{FK} Furstenberg, H.; Kesten, H.
Products of random matrices. 
Ann.\ Math.\ Statist.\ 31 (1960), 457--469.  
\bibitem{G75} Goldsheid, I.\ Ya.
Asymptotic behaviour of a product of random matrices that depend on a parameter. (Russian)
Dokl.\ Akad. Nauk SSSR 224 (1975), no. 6, 1248--1251.
\bibitem{G80a} Goldsheid, I.\ Ya.
Asymptotic properties of the product of random matrices depending on a parameter. Multicomponent random systems, pp.\ 239--283, 
Adv.\ Probab.\ Related Topics, 6, Dekker, New York, 1980. 
\bibitem{GS2}  Goldsheid, I.; Sodin, S.
Sets of  non-Lyapunov behaviour for scalar and matrix Schr\"odinger  cocycles.
Int.\ Math.\ Res.\ Not.\ IMRN 2024, no.~9, 7421--7444.
\bibitem{GorKl} Gorodetski, A.,  Kleptsyn, V.
Parametric Furstenberg theorem on random products of $\operatorname{SL} (2,\mathbb {R})$ matrices.
Adv.\ Math.\ 378 (2021), Paper No.~107522, 81 pp. 
\bibitem{Cot} Haeming, Lian. 
On the bandwidths of periodic approximations to discrete Schr\"odinger operators.
J.\ Anal.\ Math.~153 (2024), no. 2, 489–517.
\bibitem{HM} Hislop, Peter D.; M\"uller, Peter.
A lower bound for the density of states of the lattice Anderson model. 
Proc.\ Amer.\ Math.\ Soc.~136 (2008), no. 8, 2887--2893. 
\bibitem{JZ} Jitomirskaya, Svetlana; Zhu, Xiaowen.
Large deviations of the Lyapunov exponent and localization for the 1D Anderson model. 
Comm.\ Math.\ Phys. 370 (2019), no.~1, 311--324. 
\bibitem{Jeske} Jeske, F.
\"Uber lokale Positivit\"at der Zustandsdichte zuf\"alliger Schr\"odinger-Operatoren, Ph.D. thesis, Ruhr-Universit\"at Bochum, Germany, 1992
\bibitem{KlQ} Kleptsyn, V.; Quintino, F.
Phase transition of capacity for the uniform $G_\delta$-sets.
Potential Anal.\ 56 (2022), no.~4, 597--622.
\bibitem{McD} McDiarmid, Colin. On the method of bounded differences. Surveys in combinatorics, 1989 (Norwich, 1989), 148--188, 
London Math.\ Soc.\ Lecture Note Ser., 141, Cambridge Univ.\ Press, Cambridge, 1989.

\bibitem{Min} Minami, N.
Local fluctuation of the spectrum of a multidimensional Anderson tight binding model. Commun.\ Math.\ Phys.~177  (1996), 709--725.
\bibitem{Q} Quintino, F. 
Logarithmic capacity of random $G_\delta$-sets.
arXiv:2012.01593
\bibitem{Simon}
 Simon, Barry.
 Equilibrium measures and capacities in spectral theory. 
Inverse Probl.\ Imaging~1 (2007), no.~4, 713--772.
\bibitem{Weg} Wegner, F.
Bounds on the density of states in disordered systems, 
Z.\ Phys.~B 44 (1981), 9--15.
\end{thebibliography}
\end{document}